\documentclass[12pt,twoside]{amsart}
\usepackage{amsmath}
\usepackage{amssymb}
\usepackage{amscd}
\usepackage[all]{xy}
 \CompileMatrices
\title[Generalized orientations and the Bloch invariant]
{Generalized orientations and the Bloch invariant}

\author{Michel Matthey}
\thanks{}
\address{}

\author{Wolfgang Pitsch}

\author{J\'er\^ome Scherer}
\thanks{The second and third authors are supported by the program Ram\'on y Cajal, MEC,
Spain, and MEC grant MTM2004--06686. The third author was
partially supported by the Mittag-Leffler Institute in Sweden.}
\address{
\hfill\break Departament de Matem\`atiques\\
\hfill\break Universitat Aut\`onoma de Barcelona\\
\hfill\break E--08193 Bellaterra\\
\hfill\break Spain.}
\email{pitsch@mat.uab.es}
\email{jscherer@mat.uab.es}

\subjclass[2000]{Primary 57M27; Secondary 19E99, 55N20, 55P43,
57M50}

\newcommand{\im}{\operatorname{\rm Im}\nolimits}

\newcommand{\A}{\ifmmode{\cal A}\else$\cal A$\fi}
\newcommand{\Ab}{\ifmmode{{\cal A}b}\else${\cal A}b$\fi}

\newcommand{\Ss}{{\mathbb S}}

\newcommand{\C}{{\mathbb C}}
\newcommand{\Q}{{\mathbb Q}}
\newcommand{\Z}{{\mathbb Z}}

\newcommand{\R}{{\mathbb R}}

\newcommand{\Kern}{\operatorname{Ker}\nolimits}

\renewcommand{\im}{\operatorname{Im}\nolimits}

\newtheorem{theorem}{Theorem}[section]
\newtheorem{proposition}[theorem]{Proposition}
\newtheorem{corollary}[theorem]{Corollary}
\newtheorem{lemma}[theorem]{Lemma}
\newtheorem{definition}[theorem]{Definition}
\newtheorem{remark}[theorem]{Remark}
\newtheorem{example}[theorem]{Example}

\begin{document}

\begin{abstract}
For compact hyperbolic $3$-manifolds we lift the Bloch invariant
defined by Neumann and Yang to an integral class in $K_3(\C)$.
Applying the Borel and the Bloch regulators, one gets back the
volume and the Chern-Simons invariant of the manifold. We also
discuss the non-compact case, in which there appears a
$\Z/2$-ambiguity.
\end{abstract}

\maketitle

\section*{Introduction}

Suppose that $\Gamma$ is a discrete group such that the
classifying space $B\Gamma$ has a model which is a closed
orientable smooth manifold $M$ of dimension $m$. Here as usual
\emph{closed} means compact and without boundary. According to the
\emph{Borel conjecture} for $\Gamma$, the \emph{diffeomorphism
type} of $M$ should be completely determined by the
\emph{isomorphism} type of $\Gamma$. Therefore the question arises
of how much of the \emph{smooth geometry} of $M$ is encoded in the
group $\Gamma$. Similarly, recall that by the celebrated
\emph{Mostow Rigidity}, if $M$ is a closed connected orientable
hyperbolic manifold of dimension $n \geq 3$, then not only the
Borel conjecture holds for $\Gamma$, but the isometry type of $M$
is also completely determined by $\Gamma$. So, in this case, the
question refines to how the \emph{metric geometry} of $M$,
typically the \emph{hyperbolic volume} $\text{vol}(M)$ or the
\emph{Chern-Simons invariant} $\text{CS}(M)$,  can be recovered
from $\Gamma$.

Such questions have been addressed for instance by Goncharov \cite{MR1649192},
Neumann and Yang \cite{MR1663915}.
In the three dimensional case, they obtained
respectively a rational algebraic $K$-theoretical invariant, and a
\emph{Bloch invariant} in the Bloch subgroup of the scissors
congruence group  of hyperbolic  $3$-space, $\mathcal{P}(\C)$. Later, Cisneros-Molina and Jones
revisited this work in \cite{MR2006404} from a homotopical point
of view, in an attempt to lift the Bloch invariant to an integral
class in $K_3(\C)$. The later is indeed a natural candidate to contain this kind of invariant. There are two
\emph{regulators} defined on $K_3(\C)$, the Borel regulator and
the Bloch regulator. The insight of Goncharov and Neumann-Yang
tells us that their values on the invariant should give back the volume
and the Chern-Simmons invariant of the manifold.

There is one constant in all three approaches: the invariant is
obtained basically by pushing a fundamental class in ordinary
homology into $\mathcal{P}(\C)$. The main tool to relate
$\mathcal{P}(\C)$ to  $K$-theory is the Bloch-Wigner exact
sequence first published by Suslin \cite{MR1092031} and by
Dupont-Sah \cite{MR662760}. One gets directly a class in the
homology of $SL_2 \C$ by considering a Spin-structure on the
hyperbolic manifold (cf. \cite{MR2006404}). To define the
invariant in $K$-theory one has to lift this fundamental class
through a Hurewicz homomorphism and this leads to an ambiguity in
the definition. In \cite{MR1649192} this ambiguity is removed by
using rational coefficients. In our context it is more natural to
view  the Bloch-Wigner exact sequence  as a part of the long exact
sequence in stable homotopy of a cofibration. Thus instead of a
Spin-structure, which yields a $KO$-orientation \cite{MR0167985},
we are lead to consider an orientation in stable homotopy theory,
and this is provided by a stable parallelization of the
(hyperbolic) $3$-manifold.

For compact manifolds our main result is:

\medskip

{\bf Theorem~A.}
\noindent {\it Let $M$ be a closed oriented hyperbolic manifold of
dimension $3$ with fundamental group $\Gamma = \pi_1(M)$. Then, to
any stable parallelization of the tangent bundle of $M$
corresponds, in a canonical way, a $K$-theory class $\gamma(M) \in
K_3(\C)$, which depends only on the underlying Spin-structure. The
hyperbolic volume of $M$ is determined by the equality
\[
\emph{\text{bo-reg}} (\gamma(M)) =
\frac{\emph{\text{vol}}(M)}{2 \pi^2}
\]
of real numbers, where $\emph{\text{bo-reg}}~:
K_3(\C) \rightarrow \R$ is the Borel regulator
for the field of complex numbers $\C$. Furthermore, for the
Chern-Simons invariant $\emph{\text{CS}}(M)$ of $M$ we have the
congruence
\[
\mu (\gamma(M)) \equiv \frac{-\text{CS}(M) + i \cdot
\emph{\text{vol}}(M)}{2 \pi^2}  \quad (\emph{\text{mod }} \Q)
\]
of complex numbers. Here $\rho$ stands for the composite
\[
\mu~: K_3(\C) \xrightarrow{\emph{bw}} \mathcal{B}(\C)
\xrightarrow{\emph{bl-reg}} \C/\Q,
\]
where \emph{bw} is the Bloch-Wigner map for the field $\C$,
$\mathcal{B}(\C)$ is the Bloch group of $\C$, and \emph{bl-reg} is
the Bloch regulator for $\C$.}

\medskip

In the non-compact case, the problem is more intricate. The main
problem is that one has to start with a fundamental class in a
relative (generalized) homology group, and this yields naturally a
relative class. Even if we do not have to invert a Hurewicz
homomorphism we still end up with a $\Z/2$ ambiguity.

\medskip

{\bf Theorem~B.}
\noindent {\it Let $M$ be a non-compact oriented hyperbolic
manifold of dimension $3$ with finite volume. Let $\Gamma =
\pi_1(M)$ be its fundamental group. Then, to any stable
parallelization of the tangent bundle of $M$ correspond two
natural $K$-theory classes $\gamma(M)^\pm \in K_3(\C)$, which
depend only on the underlying Spin-structure. The hyperbolic
volume of $M$ is determined by the equality
\[
\emph{\text{bo-reg}} (\gamma(M)^\pm) =
\frac{\emph{\text{vol}}(M)}{2 \pi^2}.
\]
Furthermore, for the Chern-Simons invariant $\emph{\text{CS}}(M)$
of $M$ we have the congruence
\[
\rho \circ (\gamma(M)^\pm) \equiv \frac{-\text{CS}(M) + i \cdot
\emph{\text{vol}}(M)}{2 \pi^2}  \quad (\emph{\text{mod }} \Q).
\]
}

\medskip

The plan of the article is the following. Section~\ref{parallel}
is a short reminder on the theory of orientations of manifolds.
Section~\ref{sec Bloch-Wigner} is devoted to the Bloch-Wigner
exact sequence. Theorem~A collects the results of
Theorem~\ref{Bloch} and Corollary~\ref{regulator}, which are
proved in Section~\ref{sec compact}. The non-compact case,
Theorem~B, is the object of Section~\ref{sec non-compact}. Our
original plan was to construct an invariant in algebraic
$K$-theory of the group ring $\Z\Gamma$. The fact that the
Bloch-Wigner exact sequence can be reformulated in stable homotopy
simplified the construction. Due to the intimate relation of
$K\Z_*(\Z\Gamma)$ with the Isomorphism Conjectures,
\cite{MR1179537}, we decided to include our original construction
in Appendix~\ref{sec KZorientation}.

\medskip

We started this project in February 2005, but the paper was
completed only after the first author's death. It is dedicated to
the memory of our friend Michel Matthey.

\medskip

{\bf Acknowledgements}. We would like to thank Joan Porti, Jos\'e
Burgos, and Johan Dupont for enlightening discussions.

\section{Parallelizations and orientations}
\label{parallel}

Let $M$ be a closed compact connected smooth manifold of
dimension~$d$. We explain in this section the relationship between
stable parallelizations of the tangent bundle of $M$ and
orientations of $M$ with respect to the sphere spectrum
$\mathbb{S}$. For manifolds there are two ways to view
orientations. The first one, rising from orientations of vector
bundles, is cohomological in essence and the second one, rising
from patching local compatible orientations, is homological in
essence.  Both definitions agree via the so-called $S$-duality. We
call a manifold \emph{orientable}  if it is so in the classical
sense (i.e. with respect to the Eilenberg-McLane spectrum $H\Z$).
In this section $E$ denotes a ring spectrum with unit $\varepsilon
: \Ss \rightarrow E$.

\subsection{Cohomological definition}\label{subsec cohodef}
Let $\nu_M$ be the stable normal bundle of $M$ and $Th(\nu_M)$ its
Thom spectrum. For each $m \in M$, consider the map from the Thom
spectrum of this point induced by the inclusion $j_m : \Ss
\rightarrow Th(\nu_M)$. An \emph{$E$-orientation} of $M$ is a
class $t \in E^0(Th(\nu_M))$ such that for some (and hence every)
point $m \in M$ $j^{\ast}_m(t) = \pm \varepsilon \in \pi_0(E)
\cong E^0(\Ss)$.

A particularly convenient setting is when the manifold is stably
parallelizable, i.e. its normal bundle is stably trivial (and
hence its tangent bundle also). A given parallelization $\iota$
provides a trivialization of the Thom spectrum of the normal
bundle of $M$:
$$
DT(\iota): Th(\nu_M) \xrightarrow{\simeq} \Sigma^\infty M_+\, .
$$
By collapsing $M$ to a point we obtain hence a map $Th(\nu_M)
\rightarrow \Ss$ to the sphere spectrum representing a cohomology
class in $\Ss^0(Th(\nu_M))$. Composing with the unit $\varepsilon
: \Ss \rightarrow E$ we get an $E$-orientation.

\begin{example}
\label{ex 3mancoho}
{\rm Recall Stiefel's result that any orientable $3$-manifold
admits stable parallelizations (see
\cite[Problem~12-B]{MR0440554}) i.e. trivializations of the stable
tangent bundle $\tau: M \rightarrow BO$. As these correspond to
lifts of the map $\tau$ to the universal cover $EO$ up to
homotopy, one can apply obstruction theory to count them. Lifts to
the $1$-skeleton correspond to classical orientations and there
are $H^0(M; \Z/2\Z)$ possible choices. Further lifts to the
$2$-skeleton correspond to Spin-structures, and there are $H^1(M;
\Z/2\Z)$ choices at this stage. Finally, to lift further across
the $3$-skeleton one gets $H^3(M;\Z)$ choices, the so called
$p_1$-structures, where $p_1$ stands for the first Pontrjagin
class.}
\end{example}

\subsection{Homological definition}\label{subs homoldef}
A \emph{fundamental class} for $M$ with respect to
the homology theory $E$ is a class $t \in E_d(M)$ such that for
some (and therefore every) point $m \in M$ the image of $t$ in
$E_d(M,M-m) \simeq \tilde{E}_d(S^d) \simeq \tilde{E}_0(S^0) =
\pi_0(E)$ is $\pm \varepsilon$. Notice in particular that the unit $\varepsilon
: \Ss \rightarrow E$ canonically provides fundamental classes for all
spheres $S^d$.

\begin{example}
\label{stablehomotopy}
{\rm Consider the sphere spectrum $\Ss$. Then the corresponding
reduced homology theory is stable homotopy, $\tilde\Ss_n(X) \cong
\pi_n^S(X)$. An $\Ss$-orientation for $M$ is thus an element in
$\Ss_d(M)$ with the property that its image in $\Ss_d^S(M, M-m)
\cong \pi_d^S(S^d) \cong \Z$ is a generator. }
\end{example}

\subsection{$S$-duality}
\label{subsec S-duality}

We now turn to the connection between the homological and
cohomological point of view. We adopt the point of view of Rudyak
\cite{MR1627486} on $S$-duality, for another point of view see
Switzer \cite{MR1886843} or Adams \cite{MR0402720}.

\begin{definition}\label{def sduality}
{\rm Let $A, A^\ast$ be two spectra. A \emph{duality morphism} or
\emph{duality} between $A$ and $A^\ast$ is a map of spectra
$u~:\Ss \rightarrow A \wedge A^\ast$ such that for every spectrum
$E$ the following homomorphisms are isomorphisms~:
\[
\begin{array}{rcl}
u_E~: [A,E] & \longrightarrow & [\Ss, E \wedge A^\ast] \\
\phi & \longmapsto & (\phi \wedge 1_{A^\ast}) \circ u
\end{array}
\]
\[
\begin{array}{rcl}
u^E~: [A^\ast,E] & \longrightarrow & [\Ss, A \wedge E ] \\
\phi & \longmapsto & (1_A \wedge \phi ) \circ u
\end{array}
\]
The spectra $A$ and $A^\ast$ are said to be \emph{$S$-dual}. Two
spectra $A$ and $B$ are called \emph{$n$-dual}, where $n \in \Z$,
if $A$ and $\Sigma^{n} B$ are $S$-dual.}
\end{definition}

\begin{definition}\label{def nSduality}
{\rm Fixing two duality maps $u: \Ss \rightarrow A\wedge A^\ast$
and $v: \Ss \rightarrow B \wedge B^\ast$, the \emph{$S$-dual} of a
map $f: A \rightarrow B$ is then the image $f^\ast: B^\ast
\rightarrow A^\ast$ of $f$ under the isomorphism:
\[
\xymatrix{ {D: [A,B]} \ar[r]^{u_B} & {[\Ss, B \wedge A^\ast]}
\ar[r]^{{(v^{A^\ast})^{-1}}} & {[B^\ast, A^\ast].} }
\]
}
\end{definition}

In particular $f \in [A,B]$ is $S$-dual to $g \in [B^\ast,A^\ast]$
if and only if $u_B(f) = v^{A^\ast}(g)$.

\begin{example}
{\rm For any integer $n$ the spectra $\Ss^n$ and $\Ss^{-n}$ are
$S$-dual. The duality map is simply the canonical  equivalence
$\Ss \rightarrow \Ss^{n} \wedge \Ss^{-n}$.}
\end{example}

\subsection{Orientations and $S$-duality for manifolds}\label{subsec orandsforman}

For closed manifolds $S$-duality was defined by Milnor-Spanier in
\cite{MR0117750}. As we will need the precise form of the duality
map we give it in detail. Choose an embedding $M \hookrightarrow
S^N$ into a high-dimensional sphere and let $U$ be a tubular
neighborhood of $M$. The open manifold $U$ can be viewed as the
total space of the normal disc bundle of $M$, and the quotient
$\overline{U}/
\partial U$ is therefore a Thom space for the normal
bundle. Denote by $p: \overline{U} \rightarrow M$ the projection
and by $\Delta: \overline{U} \rightarrow \overline{U} \times M$
the map  $\Delta(a) = (a,p(a))$. Then $\Delta$ induces a map
$\Delta': \overline{U}/\partial U \rightarrow
\overline{U}/\partial{U} \wedge M_+$. Denote by $C: S^N
\rightarrow \overline{U}/\partial U$ the map induced by collapsing
the complement of $U$ into a point. Then we have a map $f~: S^N
\xrightarrow C \overline{U}/\partial U \xrightarrow \Delta'
(\overline{U}/\partial U) \wedge M_+$. The duality morphism is
then
$$
u = \Sigma^{-N} \Sigma^{\infty} f: \Ss \rightarrow Th \nu_M \wedge
\Sigma^{-d} \Sigma^\infty M_+.
$$
It induces the duality bijection $u_E: [Th(\nu_M) , E] \rightarrow
[\Ss, E \wedge \Sigma^{-d} \Sigma^{\infty} M_+]$ for any
spectrum~$E$.

\begin{theorem} \cite[Corollary~V.2.6]{MR1627486}
\label{thm equivcoho-homoor}
Let $M$ be a closed $E$-orientable manifold. The duality map
constructed above yields a bijective correspondence between
cohomological orientations of $M$ and fundamental classes of $M$
with respect to $E$.\hfill{\qed}
\end{theorem}

\subsection{The case of $3$-manifolds}\label{subses case3man}

In Example \ref{ex 3mancoho} we have seen that $3$-manifolds are
orientable in the cohomological sense. Therefore by Theorem
\ref{thm equivcoho-homoor} they admit fundamental classes. We
describe now the relationship between parallelizations and
homological orientations for $3$-manifolds. Since we counted the
former in Example~\ref{ex 3mancoho} we will first count the later.

\begin{lemma}
\label{AHSS}
Let $M$ be an orientable closed manifold of dimension $3$. The
Atiyah-Hirzebruch spectral sequence  for the stable homotopy of
$M$ collapses at $E^2$.
\end{lemma}
\begin{proof}
The spectral sequence is concentrated on the first four columns of
the first quadrant. The first column $H_0(M; \Ss_q) \cong \pi_q^S$
always survives to $E^\infty$ since a point is a retract of $M$.
Since $M$ is $\Ss$-orientable, the suspension spectrum of the
$3$-sphere is a retract of $\Sigma^\infty M$, so that the fourth
column $H_3(M; \Ss_q) \cong \pi_q^S$ also survives. Therefore all
differentials must be zero.
\end{proof}

\begin{proposition}
\label{numberS}
Let $M$ be an orientable closed $3$-manifold. Fundamental classes
of $M$ with respect to $\Ss$ are parametrized by
$\pi_3^s(\mathbb{S}) \times H_1(M; \Z/2\Z) \times H_2(M; \Z/2\Z)
\times \{ \pm 1 \}$.
\end{proposition}

\begin{proof}
This follows from the previous lemma since  the homomorphism
$\Ss_3(M) \rightarrow \Ss_3(M, M-m)$ can be identified with the
edge homomorphism $\Ss_3(M) \rightarrow H_3(M; \Z)$. Fixing an
orientation tells us that the image of $t$ must be a fixed
generator of $H_3(M; \Z)$.
\end{proof}

\begin{example}
\label{sphere}
{\rm There are precisely $2 \cdot |\pi_3^s(\Ss)| = 48$ different
orientations of the sphere $S^3$ with respect to stable
homotopy.}
\end{example}


If an $\Ss$-orientation of $M$ is given, a change of
trivialization can be used to modify the class in $\Ss_3(M)$ via
the Dold-Thom isomorphisms:
$$
\Ss_3(M) \xrightarrow{DT(\iota)^{-1}} \Ss_3(Th(\nu_M))
\xrightarrow{DT(\iota')} \Ss_3(M).
$$

\begin{lemma}
\label{imageJ}
Given two stable parallelizations of $S^3$ which differ only by a
$p_1$-structure $\alpha \in H^3(S^3; \Z)$, the corresponding
$\Ss$-orientations differ then by $J\alpha$, where $J: \Z \cong
\pi_3 SO \twoheadrightarrow \pi_3^S \cong \Z/24$ is the stable
$J$-homomorphism.
\end{lemma}

\begin{proof}
The change of trivialization is controlled by a map between total
spaces of trivial bundles $S^3 \times \R^N \rightarrow S^3 \times
\R^N$, for some large integer $N$. At the level of Thom spaces we
get a homotopy equivalence $f: S^{N+3} \vee S^N \rightarrow
S^{N+3} \vee S^N$. Fix the canonical $\Ss$-orientation $t$
corresponding to the inclusion $S^{N+3} \rightarrow S^{N+3} \vee
S^N$ in $\pi_{N+3}(S^{N+3} \vee S^N) \cong \pi_3^S(S^3_+) \cong
\Ss_3(S^3)$ and modify it by $f$. The edge homomorphism $e:
\Ss_3(S^3) \rightarrow \pi_3^S(S^3)$ takes both $t$ and $ft$ to
$1$, but the element in $\Kern e$ is null for $t$ and, for $ft$,
is given by the map
$$
S^{N+3} \xrightarrow{i_1} S^{N+3} \vee S^N \xrightarrow{\ f\ }
S^{N+3} \vee S^N \xrightarrow{p_2} S^N.
$$
This map is determined by its homotopy cofiber, a two cell complex
which is seen to be homotopy equivalent to $S^N \cup_{J\alpha}
e^{N+4}$, see \cite[Lemma~10.1]{MR0198470}. We conclude then since
$J$ is an epimorphism in dimension 3,
\cite[Theorem~1.5]{MR0198470}.
\end{proof}

\begin{proposition}
\label{kernel}
Let $M$ be an oriented, closed $3$-manifold. The
$\Ss$-orientations of $M$ obtained from the stable
parallelizations may differ by an arbitrary element of $\Z/24
\cong \pi_3^S \subset \Kern e$.
\end{proposition}

\begin{proof}
One obtains both stable parallelizations and $\Ss$-orientations
for $S^3$ from the ones for $M$ by collapsing the $2$-skeleton.
\end{proof}

\section{The Bloch-Wigner exact sequence}
\label{sec Bloch-Wigner}

In this section we identify the Bloch-Wigner exact sequence with
an exact sequence in stable homotopy whereas the classical point
of view is homological.

\subsection{Scissors congruence group of hyperbolic $3$-space}
\label{subsec scissors}
A standard reference fore this section is Dupont-Sah
\cite{MR662760}, see also Dupont \cite{MR1832859} or Suslin
\cite{MR1092031}. Denote by $\text{Isom}^+(\mathcal{H}^3)$ the
group of orientation-preserving isometries of the hyperbolic
$3$-space $\mathcal{H}^3$.

\begin{definition}\label{def scissors}
{\rm The \emph{scissors congruence group}
$\mathcal{P}(\mathcal{H}^3)$ is the free abelian group of symbols
$[P]$ for all polytopes $P$ in $\mathcal{H}^3$, modulo the
relations:
\begin{enumerate}
\item $[P] - [P'] -[P'']$ if $P = P' \cup P''$ and $P' \cap P''$
has no interior points;
\item $[gP] - [P]$ for $g \in \text{Isom}^+(\mathcal{H}^3)$.
\end{enumerate}}
\end{definition}

One defines  analogously $\mathcal{P}(\overline{\mathcal{H}}^3)$
where one allows some vertices of the polytopes to be ideal points
and $\mathcal{P}(\partial \mathcal{H}^3)$ where the polytopes are
all ideal polytopes (actually there is a subtlety with the later
group, see \cite[Chapter 8]{MR1832859}). Finally there is a more
algebraic description of these groups.

\begin{definition}\label{def pdc}
\rm{ Let $\mathcal{P}(\C)$ denote the abelian group generated by
symbols $z  \in \C-\{0,1\}$ and satisfying, for $z_1 \neq z_2$,
the relations:
\[
z_1 -z_2 + \frac{z_2}{z_1} - \frac{1-z_2}{1-z_1} +
\frac{1-z_2^{-1}}{1-z_1^{-1}}.
\]}
\end{definition}

The four groups are related by:

\begin{theorem}\cite[Corollary 8.18]{MR1832859}\label{them iso4grps}
There are canonical isomorphisms
\[
\mathcal{P}(\mathcal{H}^3) \cong \mathcal{P}(\overline{\mathcal{H}}^3) \cong \mathcal{P}(\partial \mathcal{H}^3) \cong \mathcal{P}(\C)^-,
\]
where $\mathcal{P}(\C)^-$ denotes the $(-1)$-eigenspace for
complex conjugation. \hfill{\qed}
\end{theorem}

\subsection{The Bloch-Wigner exact sequence}\label{subsec BW}

Recall that the group $\text{Isom}^+(\mathcal{H}^3)$ is isomorphic
to $PSL_2 \C = SL_2 \C / \{ \pm Id \}$. It acts naturally on the
boundary of hyperbolic $3$-space. Fix a point $x \in \partial
\mathcal{H}^3$ and denote by $P \subset SL_2(\C)$ the preimage of
a parabolic stabilizer of~$x$. As a group, $P$ is isomorphic to
the semi-direct product $\C \rtimes \C^\ast$, where $z \in
\C^\ast$ acts on $\C$ by multiplication by $z^2$. These groups are
all considered only as dicrete groups. Let us then denote by
$Cof(i_P)$ the homotopy cofibre of the map $i_P: BP \rightarrow
BSL_2(\C)$. The following is an integral analogue of
\cite[Lemma~2.14]{MR1649192}.

\begin{lemma}
\label{lem homoBP}
For $n \geq 1$ we have a commutative diagram where the vertical
arrows are the Hurewicz homomorphisms and the horizontal arrows
are induced by the projections $\C \rtimes \C^\ast
\twoheadrightarrow \C^\ast$:
\[
\xymatrix{
 \pi_n^s(B(\C \rtimes \C^\ast) \ar[r] \ar^- {\sim}[d] & \pi_n^s(B\C^\ast) \ar[d] \\
H_n(\C \rtimes  \C^\ast; \Z) \ar^- {\sim }[r] & H_n(B\C^\ast, \Z)
. }
\]
\end{lemma}

\begin{proof}
{}From the exact sequence of groups $1 \rightarrow \C \rightarrow
{\C \rtimes \C^\ast} \rightarrow {\C^\ast} \rightarrow 1$ we get a
fibration $B\C \rightarrow {B(\C \rtimes  \C^\ast)} \rightarrow
{B\C^\ast}$. We will prove that the Atiyah-Hirzebruch spectral
sequence for stable homotopy
\[
H_p(B\C^\ast, \pi_q^s(B\C)) \Rightarrow \pi_{p+q}^s(B(\C \rtimes
\C^\ast))
\]
collapses. Since the stable stems $\pi_n^S$ are torsion groups in
degree $n \geq 1$ and $\C$ is a rational vector space, the
Hurewicz homomorphism $\pi_n^s(B\C) \rightarrow H_n(B\C;\Z)$ is an
isomorphism. This identifies the above spectral sequence with the
ordinary homological spectral sequence. In particular the Hurewicz
map $\pi_n^s(B(\C \rtimes \C^\ast)) \rightarrow H_n(B(\C \rtimes
\C^\ast); \Z)$ is an isomorphism.

Now, $H_q(B\C; \Z) \cong \Lambda^q \C$ for any $q \geq 1$. An
element $n \in \C^\ast$ acts by multiplication by $n^2$ on $\C$
and therefore by multiplication by $n^{2q}$ on $H_q(B\C, \Z)$. The
map induced by conjugation in a group $G$ by an element $g$
together with the action of the same $g$ on a $G$-module $M$
induces the identity in homology with coefficients in $M$. As
$\C^\ast$ is abelian, in our case we have that multiplication by
$n^{2q}$ is the identity on $ H_p(B\C^\ast, H_q(B\C))$. But
multiplication by $n^{2q}-1$ is a isomorphism of the
$\C^\ast$-module $H_q(B\C;\Z)$, therefore $H_p(B\C^\ast,
H_q(B\C))= 0$ for $q \geq 1$.
\end{proof}

\begin{lemma}
\label{lem homoSL}
For $n \leq 3$, the Hurewicz homomorphism $\pi_n^S (BSL_2 \C)
\rightarrow H_n(SL_2 \C)$ is an isomorphism.
\end{lemma}

\begin{proof}
The group $SL_2 \C$ is perfect and $H_2(SL_2\C;\Z)$ is a rational
vector space \cite[Corollary~8.20]{MR1832859}. One concludes then
by an easy Atiyah-Hirzebruch spectral sequence argument.
\end{proof}

\begin{proposition}\label{Bloch-Wigner}
There is a commutative diagram with vertical isomorphisms and exact rows
\[
\xymatrix{ {\Q/\Z} \ar@{^{(}->}[r] \ar^- {}[d]&
{{\pi}_3^s(BSL_2\C)} \ar[r] \ar[d]  & {{\pi}_3^s (Cof(i_P))}
\ar[r]^{} \ar[d]& {{\pi}_2^s(BP)} \ar@{->>}[r] \ar^- {}[d] &
{{\pi}_2^s
(BSL_2\C)} \ar[d] \\
{\Q/\Z} \ar@{^{(}->}[r] & H_3(SL_2 \C; \Z) \ar[r] &{\mathcal{P}(\C)} \ar[r] &
 {\Lambda^2 (\C^\ast/ \mu_\C)} \ar@{->>}[r] &H_2(SL_2\C; \Z) }
\]
where the bottom row is the Bloch-Wigner exact sequence.
\end{proposition}
\begin{proof}
The stable Hurewicz homomorphism permits us to compare the long
exact sequences of the cofibration $BP \rightarrow BSL_2(\C)
\rightarrow Cof(i_P)$ in stable homotopy and in ordinary homology:
\[
\xymatrix{ {{\pi}_3^s(BP)} \ar[r] \ar^- {\cong}[d]&
{{\pi}_3^s(BSL_2\C)} \ar[r] \ar^- {\cong}[d]  & {{\pi}_3^s
(Cof(i_P))} \ar[r]^{} \ar[d]& {{\pi}_2^s(BP)} \ar[r] \ar^-
{\cong}[d] & {{\pi}_2^s
(BSL_2\C)} \ar^- {\cong}[d] \\
H_3(\C^\ast;\Z) \ar[r] & H_3(SL_2 \C; \Z) \ar[r] &H_3(Cof(i_P);\Z)
\ar[r] & H_2(\C^\ast;\Z) \ar[r] &H_2(SL_2\C; \Z) }
\]

The marked isomorphisms are given by Lemmas~\ref{lem homoBP}
and~\ref{lem homoSL}. It remains thus to compare the bottom exact
sequence with the Bloch-Wigner exact sequence. We have to return
to its computation by Suslin, \cite{MR1092031}.

Let $P_*$ be a projective resolution of $\Z$ over $SL_2(\C)$ and
consider the complex $C_*$ of $(n+1)$-uples of distinct points in
$\partial \mathcal{H}^3$, \cite[Chapter~2]{MR1832859}. The
naturally augmented complex $\epsilon: C_* \rightarrow \Z$ is
acyclic. Let us consider the truncated complex $ \tau C_\ast =
(\ker \epsilon \rightarrow C_0)$. The inclusion of this complex in
$C_*$ allows to compare two spectral sequences. The first one,
associated to the double complex $P_* \otimes \tau C_*$, yields a
kind of Wang sequence, which is nothing but the long exact of the
cofibration $BP \rightarrow BSL_2(\C) \rightarrow Cof(i_P)$. The
second one, associated to the double complex $P_* \otimes C_*$,
yields in low degrees the classical Bloch-Wigner sequence. In
particular we get isomorphisms $H_3(Cof(i_P);\Z) \cong
{\mathcal{P}(\C)}$ and $\im (H_3(\C^\ast;\Z) \rightarrow H_3(SL_2
\C; \Z)) \cong \Q/\Z$.
\end{proof}

\section{Lifting the Bloch invariant, the compact case}
\label{sec compact}
We construct in this section a class in $K_3(\C)$ for every
closed, compact, orientable hyperbolic $3$-manifold and show it
coincides with the Neumann-Yang \emph{Bloch invariant},
\cite{MR1663915}. In Section~\ref{parallel} we have seen that one
obtains from a stable parallelization of the normal bundle an
$\Ss$-orientation. Set $\Gamma = \pi_1 M$ and let us fix a
Spin-structure $\rho: \Gamma \rightarrow SL_2(\C)$.

\subsection{The invariant $\gamma(M)$}
\label{subsec gamma}
We start with an $\Ss$-orientation $t \in \Ss_3(B\Gamma)$ coming
from a stable parallelization that extends the Spin-structure
$\rho$ (recall from Example~\ref{ex 3mancoho} that $\rho$ provides
a trivialization of the normal bundle over the $2$-skeleton
of~$M$). Note that the reduced homology groups are canonical
direct factors of the unreduced ones for pointed spaces, so we
have a projection $\Ss_3(M) \twoheadrightarrow \widetilde \Ss_3(M)
\cong \pi_3^S(M)$, sending a given orientation $t \in \Ss_3(M)$ to
a \emph{reduced orientation class} $\tilde t$ in $\pi_3^S (M)$.

The idea is to use the structural map $\rho$ to obtain an element
in $\pi^S_3(BSL_2 \C)$. Then include $SL_2 \C$ into the infinite
special linear group $SL \C$. This defines for us an element in
$$
\pi^S_3(BSL\C) \cong \pi^S_3(BSL\C^+)\, .
$$

\begin{lemma}
\label{stabilization}
The stabilization map $\pi_3 BSL\C^+ \rightarrow \pi_3^S BSL\C^+$
is an isomorphism.
\end{lemma}

\begin{proof}
Since $BSL\C^+$ is simply connected, Freudenthal's suspension
theorem tells us that the stabilization homomorphism $\pi_3
BSL\C^+ \twoheadrightarrow \pi_3^S BSL\C^+$ is an epimorphism. The
infinite loop space $BSL\C^+$, being the universal cover of
$BGL\C^+$, gives rise to the $1$-connected spectrum $K\C\langle 1
\rangle$. The map of spectra $\Sigma^\infty BSL\C^+ \rightarrow
K\C \langle 1 \rangle$, adjunct of the identity, yields a right
inverse to the stabilization map, which must therefore be a
monomorphism.
\end{proof}

\begin{definition}
\label{Blochinvariant}
{\rm Let $M$ be a closed, compact, orientable hyperbolic
$3$-manifold with fundamental group $\Gamma$ (thus $M \simeq
B\Gamma$). Fix a Spin-structure $\rho: \Gamma \rightarrow
SL_2(\C)$ and a reduced stable orientation $t \in
\pi^S_3(B\Gamma)$ coming from a stable parallelization
extending~$\rho$. The element $\gamma(M)$ is then the image of
$\tilde t$ by the homomorphism
$$
\pi^S_3(B\Gamma) \xrightarrow{\rho_*} \pi^S_3(BSL_2 \C)
\xrightarrow{i} \pi^S_3(BSL \C) \cong \pi^S_3(BSL \C^+)
\xrightarrow{\cong} K_3 (\C)\, .
$$}
\end{definition}

\subsection{Independence from the $p_1$-structure}
\label{subsec independence}

The preceding definition apparently depends on the choice of the
orientation. We prove here that $\gamma(M)$ is completely
determined by the Spin-structure only.

\begin{lemma}
\label{dualityofcollapse}
Let $M$ be  a closed orientable manifold of dimension $d$ and
$c_{(2)}: Th(\nu_M) \rightarrow \Sigma^d \Ss$ be the map obtained
by collapsing the $2$-skeleton of~$M$. The $S$-dual map of
$c_{(2)}$ is then, up to sign, the map $i_c: \Sigma^{-N} \Ss
\rightarrow \Sigma^{-d} \Sigma^\infty M_+$ induced by the
inclusion of the center of the top-dimensional cell.
\end{lemma}

\begin{proof}
The two duality maps we consider are $u: \Ss \rightarrow Th(\nu_M)
\wedge \Sigma^{-d} \Sigma^\infty M_+$ and $v: \Ss \rightarrow
\Ss^d \wedge \Ss^{-d}$. By Definition~\ref{def sduality}, we have
to prove that the maps $(c_{(2)} \wedge 1_{\Sigma^{-d}
\Sigma^\infty M_+}) \circ u$ and $(1_{S^d} \wedge i_c) \circ v$
are homotopic, i.e. coincide in
$$
[\Ss, S^d \wedge \Sigma^{-d} \Sigma^\infty M_+] = [\Ss,
\Sigma^\infty M_+] = \pi_0^S(M_+) \cong \Z\,.
$$
The collapse map $M \rightarrow pt$ induces an isomorphism
$\pi_0^S(M_+) \rightarrow \pi^S_0(S^0)$ so we may postcompose with
this collapse map.

Let us compute the homotopy class of the map $(1_{S^d} \wedge i_c)
\circ v$
\[
\xymatrix{ \Ss \ar[r]  & \Ss^d \wedge \Ss^{-d} \ar[r]
\ar@/_1.5pc/[rr]|{Id} & \Ss^d \wedge \Sigma^{-d} M^+ \ar[r] &
\Ss^d \wedge \Sigma^{-d} S^0. }
\]

Since the duality map $v$ is an equivalence this is a generator of
$\pi_0^s(S^0) = \pi_0^s(M^+)$.

To compare it with $(C \wedge 1_{\Sigma^{-d} \Sigma^\infty M^+})
\circ u$, we turn back to the definition of the duality map $u$.
One sees that the above composite is the desuspension of the
stable map induced by the following map of spaces, where $N$
stands for a sufficiently large integer:
\[
\xymatrix{ S^{d+N} \ar[r]^-{c_{(2)}} & Th(\nu_M) \ar[r]^-{\Delta'}
& Th(\nu_M) \wedge M^+ \ar[r]^{Id \wedge C} & Th(\nu_M) \wedge S^0
\ar[r]^-{c_{(2)}} & S^{d+ N} \wedge S^0 \simeq S^{d+N} }
\]
This map is equal to the map induced by the collapse of complement
of the tubular neighborhood of $M$ restricted to the
top-dimensional cell. The tubular neighborhood restricted to the
$n^{th}$ cell is a trivial disc bundle, therefore the collapse map
$S^{d+ N} \rightarrow D^d \times D^N /\partial(D^d \times D^N) =
S^{d+N}$ is of degree $\pm 1$.
\end{proof}

\begin{proposition}
\label{independent}
Let $M$ be a closed, compact, orientable hyperbolic $3$-manifold.
The reduced orientation class $\tilde t \in \pi^S_3(B\Gamma)$ is
independent of the $p_1$-structure. Consequently the element
$\gamma(M)$ depends only on the Spin-structure.
\end{proposition}

\begin{proof}
We have a cofibre sequence  $Th(\nu_M \vert_{M^{(2)}}) \rightarrow
Th(\nu_M) \xrightarrow{c_{(2)}} \Sigma^3 \Ss$ of spectra of finite
type. Therefore, by \cite[Lemma II.2.10]{MR1627486}, we have an
$S$-dual cofibre sequence of spectra of finite type $\Sigma^{-3}
\Ss \rightarrow \Sigma^{-3} \Sigma^\infty M_+ \rightarrow
(Th(\nu_M) \vert_{M^{(2)}})^\ast$, where the first map has been
identified in Lemma~\ref{dualityofcollapse}.

As a consequence we have a commutative diagram, where the vertical
arrows are induced by $S$-duality :
\[
\xymatrix{ {\Ss^0 (S^3)} \ar[r] \ar[d] & {\Ss^0 (Th \nu_M)} \ar[r]
\ar[d]
& {\Ss^0 ( Th \nu_m \vert_{M^{(2)}})} \ar[d] \\
{\pi^S_3} \ar[r] & {\pi^S_3 (M_+)} \ar[r]  & {\Ss_0 (Th(\nu_M
\vert_{M^{(2)}})^\ast)} }
\]
The map $\Ss \rightarrow \Sigma^\infty M_+ $ splits so that the
bottom row is a short exact sequence and we can identify $\Ss_0
(Th(\nu_M \vert_{M^{(2)}})^\ast)$ with $\pi_3^S (M)$. The diagram
shows that the reduced orientation class $\tilde t \in \pi_3^S
(M)$ is $S$-dual to the cohomological orientation class restricted
to the $2$-skeleton, which is unaffected by a change of
$p_1$-structure.
\end{proof}

\subsection{Comparison with the Bloch invariant}
\label{subsec Bloch}

Let us recall how Neumann and Yang construct in~\cite{MR1663915}
the Bloch invariant $\beta(M) \in \mathcal B (\C)$. The later is
the kernel of the morphism $\mathcal P (\C) \rightarrow \Lambda^2
(\C^\ast/ \mu_\C)$ in the Bloch-Wigner exact sequence,
Proposition~\ref{Bloch-Wigner}. Since $M$ is oriented hyperbolic,
$\Gamma \subset PSL_2 \C$ and $M$ can be identified with the
quotient $\mathcal H^3/\Gamma$. Choose a fundamental polytope $P
\subset \mathcal H^3$ for the action of $\Gamma$ and define
$\beta(M) = [P] \in \mathcal P (\mathcal H^3)$. One can check that
$\beta(M)$ coincides with the image of the fundamental class
through the composite
$$
H_3(M; \Z) \longrightarrow H_3(PSL_2 \C; \Z) \longrightarrow
\mathcal P(\C).
$$
This proves that $\beta(M)$ is well-defined, and lies indeed in
$\mathcal B (\C)$.

\begin{theorem}
\label{Bloch}
Let $M$ be a closed, compact, orientable hyperbolic $3$-manifold.
The element $\gamma(M)$ lifts the Bloch invariant $\beta(M)$.
\end{theorem}

\begin{proof}
It is well-known that the cokernel of the natural map $K_3^M (\C)
\rightarrow K_3 (\C)$ provides a splitting for $H_3(SL_2 \C; \Z)
\rightarrow K_3 (\C)$. Moreover the morphism $H_3(SL_2 \C; \Z)
\rightarrow \mathcal P(\C)$ factors through $H_3(PSL_2 \C; \Z)$.
Therefore we have a commutative diagram
\[
\xymatrix{ \pi_3^S (M) \ar[r]^-{\rho} \ar[d] & \pi_3^S (BSL_2 \C)
\ar[r] \ar[d] & K_3(\C)  \ar[d] \\
H_3(M; \Z) \ar[r] & H_3(PSL_2 \C; \Z) \ar[r] & {\mathcal P}(\C)}
\]
and obviously the reduced $\Ss$-orientation $\tilde t$ maps to an
orientation in $H_3(M; \Z)$.
\end{proof}

\begin{remark}
\label{compare}
{\rm Our approach can be applied in higher dimensions, since the
same definition can be used in a straightforward manner to define
a class in $K_n (\C)$ associated to an $n$-dimensional
$\Ss$-oriented hyperbolic manifold. This definition might of
course depend on the chosen orientation in general, if it exists.}
\end{remark}

Borel defined in \cite{MR0506168} the Borel regulator $\hbox{\rm
bo-reg}_\C: K_3(\C) \rightarrow \R$. Likewise the Bloch regulator
is a map $\hbox{\rm bl-reg}_\C: \mathcal{B}(\C) \rightarrow \C/\Q$
and the Bloch-Wigner map is a map $\hbox{\rm bw}_\C: K_3(\C)
\rightarrow \mathcal{B}(\C)$.

\begin{corollary}
\label{regulator}
Let $M$ be a closed compact oriented hyperbolic manifold of
dimension~$3$ with fundamental group $\Gamma$. Then, to a
Spin-structure $\rho$ corresponds, in a canonical way, a class
$\gamma(M) \in K_3(\C)$ such that the hyperbolic volume of $M$ is
determined by the equality
$$
\hbox{\rm bo-reg}_\C(\gamma(M)) = \frac{\hbox{\rm
vol}(M)}{2\pi^2}.
$$
Furthermore the Chern-Simmons invariant $\hbox{\rm CS}(M)$ is
determined by the congruence
\[
\mu (\gamma(M)) \equiv \frac{-\text{CS}(M) + i \cdot
\emph{\text{vol}}(M)}{2 \pi^2}  \quad (\emph{\text{mod }} \Q).
\]
\end{corollary}

\begin{proof}
This follows directly from Theorem~\ref{Bloch}. Neumann and Yang
prove in \cite[Theorem~1.3]{MR1663915} that one can recover the
volume and the Chern-Simmons invariant via the Borel and Bloch
regulators.
\end{proof}

\section{Lifting the Bloch invariant, the non-compact case}
\label{sec non-compact}
Let $M$ be a non-compact, orientable, hyperbolic $3$-manifold of
finite volume with $\Gamma = \pi_1(M)$. Since $M$ has finite
volume it has a finite number of cusps and all of them are
toroidal, \cite[Theorem~10.2.1]{MR1299730}. Choose such a cusp $x
\in M$ and denote by $P \subset SL_2(\C)$ the preimage of the
parabolic stabilizer of~$x$. As in Subsection~\ref{subsec BW},
$i_P$ denotes the map $BP \rightarrow BSL_2 \C$. Choose a
Spin-structure on $M$, i.e. a homomorphism $\Gamma \rightarrow
SL_2(\C)$. The representation $\rho$ contains parabolic elements,
i.e. elements fixing a point in the boundary $\partial
\overline{\mathcal{H}}^3$. Choose a sufficiently small
$\delta$-horosphere around each cusp of $M$ and denote by
$M_\delta$ the compact submanifold obtained by removing these
horospheres from the cups of $M$, \cite[Theorem~4.5.7]{MR1435975}.
The inclusion $M_\delta \hookrightarrow M$ is a homotopy
equivalence.

\subsection{A first indeterminacy for $\gamma(M)$}
\label{subsec first}
Let $T \subset \partial M_\delta$ denote any component of the
boundary, so $T \simeq S^1 \times S^1$. Consider the composite
\[
\xymatrix{ T \ar@{^{(}->}[r] & \partial M_\delta \ar@{^{(}->}[r] &
M_\delta \ar^{B\rho \ \ \ \ }[r] & BSL_2(\C) \ar[r] & Cof(i_P) }
\]
As the action of $SL_2(\C)$ is transitive on the boundary of the
hyperbolic space, all stabilizers of points in $\partial
\mathcal{H}^3$ are conjugate. The inclusion of $\Z \oplus \Z$ into
$SL_2(\C)$ is then conjugate to an inclusion into $P$, so that the
map $T \rightarrow Cof(i_P)$ is null-homotopic.

So, from the choice of the Spin-structure, we get a map
$M_{\delta}/\partial M_{\delta} \rightarrow Cof(i_P)$, which is
well-defined up to homotopy. A stable parallelization of the
tangent bundle of $M_\delta$ gives rise to a fundamental class $t
\in \Ss_3(M_\delta, \partial M_\delta) \cong
\pi_3^S(M_\delta/\partial M_\delta)$. Pushing this class by the
above map, we get a well-defined class $\gamma_P(M) \in
\pi^S_3(Cof(i_P))$.

\begin{theorem}
\label{indeterminacy}
Let $M$ be a non-compact, orientable, hyperbolic $3$-manifold of
finite volume. It is then always possible to lift the class
$\gamma_P(M)$ to a class $\gamma(M) \in K_3(\C)$, and there are
$\Q/\Z$ possible lifts.
\end{theorem}

\begin{proof}
According to Proposition~\ref{Bloch-Wigner}, the class
$\gamma_P(M)$ lives in ${\mathcal{P}(\C)} \cong
\pi^S_3(Cof(i_P))$. Thus our invariant $\gamma_P(M)$ coincides in
fact with the Bloch invariant $\beta(M)$, defined in an analgous
way to the compact case. We wish to lift it through the connecting
morphism $\delta: \pi_3^S BSL_2(\C) \rightarrow \pi_3^S
(Cof(i_P))$.

According to Neumann and Yang, \cite[Section~5]{MR1663915} the
Bloch invariant is the scissors congruence class of any hyperbolic
ideal triangulation of $M$ and this class belongs to the kernel
$\mathcal B (\C)$ of $\mathcal{P}(\C)) \rightarrow \Lambda^2
(\C^\ast/ \mu_\C)$.

The existence of the lift follows at once from
Proposition~\ref{Bloch-Wigner}. This explains the $\Q/\Z$
indeterminacy: the image of the map $\pi_3^S(BP) \rightarrow
\pi_3^S(BSL_2(\C))$ is isomorphic to $\Q/\Z$.  Now it suffices to
push any lift to $\pi_3^S(BSL(\C))$, a group isomorphic to
$K_3(\C)$ by Lemma~\ref{stabilization}.
\end{proof}

\begin{remark}
\label{rem geominter}
{\rm The fact that the Bloch invariant lies in $\mathcal B (\C)$
has a nice geometrical interpretation. Hyperbolic tetrahedra up to
isometry are in one to one correspondance with elemnts of $\C
-\{0,1\}$ , the modulus of the tetrahedron. If one starts with a
collection of such tetrahedra and wants to glue them to a
hyperbolic space then a theorem of Thurston says that the moduli
of the tetrahedra have to satisfy a compatibility relation in
$\Lambda^2 (\C-\{0,1\})$, namely $\Sigma (z\wedge(1-z)) = 0$. The
above morphism $\mathcal{P}(\C) \rightarrow \Lambda^2 (\C^\ast/
\mu_\C)$ is $z \mapsto 2(z \wedge (1-z))$. In particular the image
under this morphism of an ideal triangulation of the hyperbolic
manifold $M$ will be trivially~$0$ since we started with an
hyperbolic manifold.}
\end{remark}

Theorem~\ref{indeterminacy} immediately provides the following.

\begin{corollary}\cite[Theorem~1.1]{MR1649192}
\label{goncharov}
Let $M$ be a non-compact, orientable, hyperbolic $3$-manifold of
finite volume. Then $M$ defines naturally a class $\gamma(M) \in
K_3(\C) \otimes \Q$ such that $bo(\gamma(M)) = vol(M)$.
\hfill{\qed}
\end{corollary}

\subsection{Reducing the indeterminacy}
\label{subsec reduce}

To reduce the $\Q/\Z$ indeterminacy in the non-compact case one
can make use of the following fact. In the compact case we  do get
a class in $\pi_3^s(BSL_2 \C)$ which is further stabilized to
$\pi_3^s(BSL\C)) \simeq K_3(\C)$. Denote by $\tau$ the involution
of $\pi_3^s(BSL_2\C)$ induced by complex conjugation.

\begin{proposition}\label{prop decok3}
If $\pi_3^s(BSL_2\C)^\pm = \ker (1 \mp \tau)$ then
\begin{enumerate}
\item $\pi_3^s(BSL_2\C) = \pi_3^s(BSL_2\C)^+  + \ \pi_3^s(BSL_2\C)^-$;
\item $\pi_3^s(BSL_2\C)^+ \cap \pi_3^s(BSL_2\C)^- \cong \Z/2$;
\item the image of $\pi_3^s(BP) \rightarrow \pi_3^s(BSL_2\C)$ lies
in $\pi_3^s(BSL_2\C)^+$.
\end{enumerate}
\end{proposition}
\begin{proof}
By Lemma~\ref{lem homoSL} $\pi_3^s(BSL_2 \C) \cong H_3(SL_2 \C;
\Z)$, so it is enough to prove the assertion in homology.
According to \cite[Corollary~8.20]{MR1832859} the group $H_3(SL_2
\C; \Z)$ is divisible, so that any element $c$ can be written
$\frac{c+\tau(c)}{2} + \frac{c - \tau(c)}{2} \in H_3(SL_2\C;\Z)^+
+ \ H_3(SL_2\C; \Z)^-$. This proves point~(1).

Any element in the intersection $H_3(SL_2\C;\Z)^+  \cap
H_3(SL_2\C; \Z)^-$ is $2$-torsion. The computations in
\cite[Corollary~8.20]{MR1832859} show that the torsion in
$H_3(SL_2\C; \Z)$ is isomorphic to $\Q/\Z$. Point~(2) follows.

The torsion subgroup of $H_3(SL_2\C; \Z)$ is the image of the
composite $H_3(S^1; \Z) \rightarrow H_3(P; \Z) \rightarrow
H_3(SL_2\C; \Z)$. The action of $\tau$ on the subgroup $S^1$ in
$P$ coincides with the action induced by conjugation in $SL_2(\C)$
by the matrix $\left(
\begin{matrix} 0 & -1 \\ 1 & 0 \end{matrix} \right)$. Point (3)
follows since conjugation by an element of a group
induces the identity.
\end{proof}

For any compact hyperbolic manifold its invariant $\gamma(M)$ lies
in $\pi_3^s(BSL_2\C)^-$. Indeed the following diagram commutes:
\[
\xymatrix{H_3(M; \Z) \ar[r]^{\rho_*} \ar[d]_{-1} & H_3(SL_2\C; \Z)
\ar[d]^{\tau_*} \\
H_3(M; \Z) \ar[r]^{\rho_*} & H_3(SL_2\C; \Z) }
\]

In view of point $(3)$ in the above proposition it is natural to
choose as lifting a class in $\pi_3^s(BSL_2\C)^-$, and this
reduces the ambiguity to $\Z/2$.

\begin{theorem}
\label{thm indeterminacy}
Let $M$ be a non-compact, orientable, hyperbolic $3$-manifold of
finite volume. There are two natural lifts $\gamma(M)^\pm \in
K_3(\C)$ of the class $\gamma_P(M)$. \hfill{\qed}
\end{theorem}

\appendix
\section{Orientation with respect to algebraic $K$-theory}
\label{sec KZorientation}
To generalize this approach to higher dimensional manifolds, one
cannot follow the same strategy, as it is not known whether or not
all hyperbolic manifolds are stably parallelizable. There is
however an intermediate way, between stable homotopy and ordinary
homology. What we have done in the three dimensional situation was
to start with an $\Ss$-orientation, and the former approaches
\cite{MR1649192}, \cite{MR1663915}, and \cite{MR2006404} all
roughly started from the fundamental class in homology.

The first author's original insight to the question of lifting the
Bloch invariant was to work with $K\Z$-orientation, where $K\Z$
denote the connective spectrum of the algebraic $K$-theory of the
integers. We believe that this is a point of view which is close
enough to ordinary homology (or topological $K$-theory) so as to
be able to do computations, but at the same time not too far away
from the stable homotopy so that the above techniques to construct
an invariant in $K_3(\C)$ can go through.

\medskip

In his foundational paper \cite{MR0447373} Loday defines a product
in algebraic $K$-theory by means of a pairing of spectra (in the
sense of Whitehead). Given two rings $R$ and $S$, consider the
connective $\Omega$-spectra $KR$ and $KS$ corresponding to the
infinite loop spaces $BGLR^+ \times K_0 R$ and $BGLS^+ \times K_0
S$ respectively (the deloopings are given by the spaces
$BGL(S^nR)^+$ where $SR$ denotes the suspension of the ring $R$).
Then there exists a pairing
$$
\star: KS \wedge KR \rightarrow K(S \otimes R).
$$
We will be interested in the case when $S = \Z$. In this case the
pairing goes to $KR$. The pairing includes in particular
compatible maps
$$
BGL(S^n\Z)^+ \wedge BGLR^+ \rightarrow BGL(S^n\Z \otimes R)^+ =
BGL(S^n R)^+
$$
which yield a map of spectra $\star: K\Z \wedge BGLR^+ \rightarrow
KR$. In order to compare the present construction with the
previous one based on an $\Ss$-orientation, we will need to
understand the map obtained by precomposing with $\varepsilon
\wedge 1$, where $\varepsilon: \Ss \rightarrow K\Z$ is the unit of
the ring spectrum $K\Z$. We first look at the global pairing of
spectra.

\begin{lemma}
\label{May}
The composite map $\Ss \wedge KR \xrightarrow{\varepsilon \wedge
1} K\Z \wedge KR \xrightarrow{\star} KR$ is the identity.
\end{lemma}

\begin{proof}
We learn from May, \cite{MR593258}, that $KR$ is a ring spectrum.
In particular the composite $\Ss \wedge KR
\xrightarrow{\varepsilon \wedge 1} KR \wedge KR
\xrightarrow{\star} K(R \otimes R) \xrightarrow{\mu} KR$ is the
identity. By naturality and using the inclusion $\Z
\hookrightarrow \C$ we see that the map from the statement must be
the identity as well.
\end{proof}

We are interested in the infinite loop space $BGLR^+$ and wish to
compare it to the spectrum $KR$. For that purpose we use the pair
of adjoint functors $\Sigma^\infty: Spaces \leftrightarrows
Spectra: \Omega^\infty$, where $\Sigma^\infty X = \Ss \wedge X$ is
the suspension spectrum of the space $X$ and $\Omega^\infty E$ is
the $0$th term of the $\Omega$-spectrum representing the
cohomology theory $E^*$. If $E$ is an $\Omega$-spectrum, then
$\Omega^\infty E = E_0$ and we write $a:\Ss \wedge E_0 \rightarrow
E$ for the adjoint of the identity.

\begin{proposition}
\label{triangle}
The composite map $\Ss \wedge BGLR^+ \xrightarrow{\varepsilon
\wedge 1} K\Z \wedge BGLR^+ \xrightarrow{\star} KR$ is homotopic
to $a: \Ss \wedge BGLR^+ \rightarrow KR$.
\end{proposition}

\begin{proof}
We consider the commutative diagram
\[
\xymatrix{ \Ss \wedge \Ss \wedge BGLR^+ \ar[r]^{\varepsilon \wedge
1 \wedge 1} \ar[d]_{1 \wedge a} & K\Z \wedge \Ss \wedge BGLR^+
\ar[d]^{1
\wedge a} \ar[dr]^{\star} \\
\Ss \wedge KR \ar[r]_{\varepsilon \wedge 1} & K\Z \wedge KR
\ar[r]_{\star} & KR }
\]
The square is obviously commutative and the triangle commutes up
to homotopy since the Loday product $\star$ forms a Whitehead
pairing, \cite[p.346]{MR0447373}.
\end{proof}

Thus we can recover the invariant $\gamma(M)$ as follows. Consider
the composite
$$
h: K\Z \wedge M \xrightarrow{1 \wedge B\rho} K\Z \wedge BSL_2(\C)
\longrightarrow K\Z \wedge BGL(\C)^+ \xrightarrow{\star} K\C.
$$

\begin{proposition}
\label{KZgamma}
Let $M$ be a closed, compact, orientable hyperbolic $3$-manifold
and choose a $K\Z$-orientation $s \in K\Z_3(M) \cong \pi_3(K\Z
\wedge M)$. The invariant $\gamma(M) \in K_3(\C)$ is then equal to
$h_*(s)$. \hfill{\qed}
\end{proposition}

Between the $K\Z$-orientation and the invariant $\gamma(M)$ there
is an interesting class in $K_3(\Z \Gamma)$. It is obtained as the
image of the $K\Z$-orientation under the composite
$$
K\Z_3(B\Gamma) \longrightarrow K\Z_3(BGL(\Z \Gamma)^+)
\longrightarrow K_3(\Z \Gamma),
$$
where the first arrow is induced by the canonical inclusion
$\Gamma \hookrightarrow GL_1(\Z\Gamma)$ and the second is a Loday
product. It is not difficult to see that we recover $\gamma(M)$ by
further composing with
$$
K_3(\Z \Gamma) \xrightarrow{\rho_*} K_3(\Z SL_2 \C)
\longrightarrow  K_3(M_2 \C) \cong K_3(\C).
$$
The second arrow is the fusion map, which takes the formal sum of
invertible matrices to the actual sum in $M_2 \C$. The final
isomorphism is just Morita invariance.



\end{document}